\documentclass[12pt]{article}
\usepackage[expansion=false]{microtype}
\usepackage[reqno]{amsmath}
\usepackage{fixmath}
\usepackage{amssymb, amsthm, enumerate}
\usepackage{mathtools}
\usepackage{graphicx}

\usepackage{hyperref}

\DeclarePairedDelimiter{\norm}{\lVert}{\rVert}

\newtheoremstyle{plainsl}%
    {\topsep}
    {\topsep}
    {\slshape} 
    {}
    {\normalfont\bfseries}
    {.}
    { }
    {}

\swapnumbers

{\theoremstyle{plainsl}
\newtheorem{theorem}{Theorem}[section]
\newtheorem{lemma}[theorem]{Lemma}
\newtheorem{corollary}[theorem]{Corollary}}
{\theoremstyle{remark}
}

\newcommand\lref[1]{Lemma~\ref{lem:#1}}

\renewcommand\proof{\noindent\textsl{Proof. }}
\newcommand\sqr[2]{{\vbox{\hrule height.#2pt
    \hbox{\vrule width.#2pt height#1pt \kern#1pt
        \vrule width.#2pt}\hrule height.#2pt}}}
\renewcommand\qed{%
    \ifmmode\eqno\sqr53
    \else\nolinebreak\ \hfill\sqr53\medbreak\fi}

\numberwithin{equation}{section}

\newcommand\al{\alpha}
\newcommand\be{\beta}

\newcommand\ga{\gamma}

\renewcommand\th{\theta} 

\newcommand\cx{{\mathbb C}}

\newcommand\ints{{\mathbb Z}}
\newcommand\re{{\mathbb R}}
\newcommand\rats{{\mathbb Q}}
\newcommand\comp[1]{{\mkern2mu\overline{\mkern-2mu#1}}}
\newcommand\diff{\mathbin{\mkern-1.5mu\setminus\mkern-1.5mu}}

\newcommand\seq[3]{#1_{#2},\ldots,#1_{#3}}

\title{Pretty good state transfer on double stars}

\author{Xiaoxia Fan\thanks{
School of Mathematics and Statistics,
 Lanzhou University,  Lanzhou, Gansu, 730000, China.} \thanks{Currently a visiting Ph.D. student at the
 Department of Combinatorics and Optimization in University of Waterloo from September 2010 to September 2012; partially supported by NSF of China (No.10971086).}
\ \  Chris Godsil\thanks{
    Combinatorics \& Optimization
    University of Waterloo, N2L 3G1, Waterloo, Canada.}\\
    \texttt{x25fan@gmail.com; cgodsil@uwaterloo.ca}}

\begin{document}
\maketitle

\begin{abstract}
	Let $A$ be the adjacency matrix of a graph $X$ and suppose $U(t)=\exp(itA)$. We view
	$A$ as acting on $\cx^{V(X)}$ and take the standard basis of this space to be the
	vectors $e_u$ for $u$ in $V(X)$. Physicists say that we have
	\textsl{perfect state transfer}
	from vertex $u$ to $v$ at time $\tau$ if there is a scalar $\gamma$ such that
	\[
		U(\tau)e_u = \gamma e_v.
	\]
	(Since $U(t)$ is unitary, $\norm\gamma=1$.) For example, if $X$ is the $d$-cube and
	$u$ and $v$ are at distance $d$ then we have perfect state transfer from $u$ to $v$
	at time $\pi/2$. Despite the existence of this nice family, it has become
	clear that perfect state transfer is rare. Hence we consider a relaxation: we say
	that we have \textsl{pretty good} state transfer from $u$ to $v$ if there is
	a complex
	number $\gamma$ and, for each positive real $\epsilon$ there is a time $t$ such that
	\[
		\norm{U(t)e_u - \gamma e_v} < \epsilon.
	\]
	Again we necessarily have $|\gamma|=1$.
	
	In Godsil, Kirkland, Severini and Smith \cite{gkss} it is shown
	that we have have pretty good state transfer between the end vertices of the
	path $P_n$ if and only $n+1$ is a power of two, a prime, or twice a prime.
	(There is perfect state transfer between the end vertices only for
	$P_2$ and $P_3$.) It is something of a surprise that the occurrence of
	pretty good state transfer is characterized by a number-theoretic condition.
	                                                                             
	In this paper we extend the theory of pretty good state transfer. We provide
	what is only the second family of graphs where pretty good state transfer occurs.
	The graphs we use are the \textsl{double-star} graphs $S_{k,\ell}$, these are trees
	with a vertex of degree $k+1$ adjacent to a vertex of degree $\ell+1$, and all
	other vertices of degree one. We prove that perfect state transfer does not occur
	in any graph in this family. We show that if $\ell>2$, then there is pretty
	good state transfer in $S_{2,\ell}$ between the two end vertices adjacent
	to the vertex of degree three. If $k,\ell>2$, we prove that there is never
	perfect state transfer between the two vertices of degree at least three, 
	and we show that there is pretty
	good state transfer between them if and only these vertices both have degree
	$k+1$ and $4k+1$ is not a perfect square. Thus we find again the the existence of
	perfect state transfer depends on a number theoretic condition. It is also 
	interesting that although no double stars have perfect state transfer, there
	are some that admit pretty good state transfer. 
\end{abstract}

\section{Introduction}

Let $X$ be a graph on $n$ vertices with adjacency matrix $A$
and let $U(t)$ denote the matrix-valued function $\exp(iAt)$. We note that
$U(t)$ is both symmetric and unitary, and that it
determines what is called a \textsl{continuous quantum walk}. Work in quantum
computing has raised many questions about the relation between physically
interesting properties of $U(t)$ and properties of the graph $X$. For recent
surveys see \cite{cg-pst}, \cite{Kay}.

The physical properties of interest to us in this paper are
\textsl{perfect state transfer} and \textsl{pretty good state transfer}.
Assume $n=|V(X)|$ and identify the coordinates of $\re^n$ with $V(X)$. If $u\in V(X)$,
we use $e_u$ to denote the standard basis vector indexed by $u$. If $u$ and $v$ are distinct
vertices of $X$ we say we have \textsl{perfect state transfer} from $u$ to $v$
at time $\tau$ if there is a complex number $\gamma$ such that
\[
	U(\tau)e_u = \gamma e_v.
\]
Since $U(t)$ is unitary, $|\gamma|=1$. The evidence is that perfect state transfer
is uncommon, and we consider a relaxation of it. We say that we have
\textsl{pretty good state transfer} from $u$ to $v$ if there is a complex number $\gamma$ and,
for each positive real $\epsilon$ there is a time $t$ such that
\[
	\norm{U(t)e_u - \gamma e_v} < \epsilon.
\]
Again we necessarily have $|\gamma|=1$. Pretty good state transfer was introduced
in \cite{cg-pst}.

There is a considerable literature on perfect state transfer. In a seminal paper
on the topic Christandl et al.~\cite{cdel-pst} show that there is perfect state transfer between
the end vertices of the paths $P_2$ and $P_3$, but perfect state transfer does not occur
between the end vertices of any path on four or more vertices. From \cite{cg-pst}
we know that, for any integer $k$, there are only finitely many connected graphs
with maximum valency at most $k$ on which perfect state transfer occurs.

Much less is known about pretty good state transfer. In \cite{cg-pst} it is shown that
it takes place on $P_4$ and $P_5$. Vinet and Zhedanov \cite{VinZh} have studied pretty good
state transfer on weighted paths with loops
Godsil, Kirkman, Severini and Smith that pretty good state transfer occurs between
the end-vertices of $P_n$ if and only $n+1=2^m$, or if
$n+1=p$ or $2p$ where $p$ is an odd prime. It is surprising to see that the existence
of pretty good state transfer depends so delicately on the prime divisors of $n+1$.

In this paper we provide a second class of graphs where pretty good state transfer
occurs if and only if a number theoretic condition holds.
Let $S_{k,k}$ denote the graph we get by taking two copies of $K_{1,k}$ and joining the
two vertices of degree $k$ by a new edge. We show that there is pretty good state transfer
between the vertices of degree $k+1$ in $S_{k,k}$ if and only if $4k+1$ is not a perfect
square. We also show that there is never perfect state transfer between these two vertices.
We conclude the paper with some remarks that show that if pretty good state transfer
does occur on a graph, then, in a sense, it must occur regularly.

\section{Quotients}

We introduce a useful tool. A more expansive treatment will be found in \cite[Ch.~9]{cggrbk}.

Let $X$ be a graph. A partition $\pi$ of vertex set $V(X)$ with
cells
\[
	C_1,C_2,\cdots,C_r
\]
is \textsl{equitable}, if the number of neighbors
in $C_j$ of any vertex $u$ in $C_i$ is a constant $b_{ij}$. The directed graph with the $r$
cells of $\pi$ as its vertices and $b_{ij}$ arcs from the $i$-th to
the $j$-th cells of $\pi$ is called the \textsl{quotient} of $X$ over
$\pi$, and denoted by $X/\pi$. The entries of the adjacency matrix
of this quotient graph are given by $A(X/\pi)_{i,j}=b_{ij}$. We can
symmetrize $A(X/\pi)$ to $B$ by letting
$B_{i,j}=\sqrt{b_{ij}b_{ji}}$. We call the (weighted) graph with adjacency
matrix $B$ the \textsl{symmetrized quotient graph}. In the following we always use $B$ to
denote the symmetrized form of the matrix $A(X/\pi)$.

If $\pi$ is a partition of $V(X)$, its characteristic matrix, denoted by $P$, is the 01-matrix
whose columns are the characteristic vectors of the cells of $\pi$, viewed as subsets of $V(X)$.
If we normalize the characteristic matrix $P$ such that each column have length one, then we
obtain the \textsl{normalized  characteristic matrix} of $\pi$, denoted by $Q$. Note that $Q^TQ=I$
and $QQ^T$ is a block diagonal matrix with diagonal blocks $\frac{1}{r}J_r$, where $J_r$ is the
all-ones matrix of order $r\times r$, and the size of the $i$-th block is the size of the $i$-th
cell of $\pi$. The vertex $u$ forms a singleton cell of $\pi$ if and only if $QQ^Te_u=e_u$.

For the sake of convenience, in the following text, we always denote by $\{u\}$ the singleton cell
$\{u\}$.

Since $A$ is symmetric, it has a spectral decomposition
\[
	A = \sum_r \theta_r E_r
\]
where $\theta_r$ runs over the distinct eigenvalues $\theta_r$ of $A$ and $E_r$ is the matrix
that represents orthogonal projection onto the the eigenspace belonging to $\theta_r$.

\begin{lemma}\label{lem:idemquo}
	Let $\pi$ be an equitable partition of $X$ with normalized characteristic matrix $Q$.
	Let $A$ be the adjacency matrix of $X$ and let $B$ be the adjacency matrix of the
	symmetrized quotient graph. Then the idempotents in the spectral decomposition of
	$B$ are the non-zero matrices $Q^TE_rQ$, where $E_r$ runs over the idempotents in
	the spectral decomposition of $A$.
\end{lemma}

\proof
As $\pi$ is equitable, $AQ=QB$. Hence $A^kQ=QB^k$ and so if $f(t)$ is a polynomial
then $f(A)Q=Q f(B)$. There is a polynomial $f_r$ such that $f_r(A)=E_r$, and hence
\[
	E_r Q = Q f_r(B).
\]
Then $f_r(B)=Q^TE_rQ$ is symmetric, we show it is idempotent. We have
\[
	(Q^TE_rQ)^2 = Q^TE_r QQ^T E_r Q
\]
and since $QQ^T$ commutes with $A$, it commutes with $E_r$. Since $Q^TQ=I$ we then have
\[
	Q^TE_r QQ^T E_r Q = Q^TQQ^T E_rQ = Q^TE_rQ.
\]
It follows that
\begin{equation}
	\label{eq:bqtaq}
	B = Q^TAQ = \sum_r \theta_r Q^TE_rQ
\end{equation}
and since
\[
	\sum_r Q^TE_rQ = Q^TQ = I
\]
we conclude that \eqref{eq:bqtaq} is the spectral decomposition of $B$.\qed

If $\{u\}$ and $\{v\}$ are singletons in $\pi$, the $2\times2$ submatrix of $Q^TE_rQ$ indexed
by $\{u\}$ and $\{v\}$ is equal to the $2\times 2$ submatrix of $E_r$ indexed by $u$ and $v$.
Since $u$ and $v$ are strongly cospectral if and only if for each idempotent $E_r$,
we have $(E_r)_{u,u} = (E_r)_{v,v} = \pm(E_r)_{u,v}$.

In \cite{YG}, R.~Bachman et el.~studied perfect state transfer of quantum walks on quotient graphs.
We state their result without proof. For the details, see \cite[Theorem~2]{YG}.

\begin{lemma}
\label{lem:gy}
	Let $X$ be a graph with an equitable partition $\pi$ and assume
	$\{a\}$ and $\{b\}$ are singleton cells of $\pi$. Let $B$ denote the adjacency matrix
	of the symmetrized quotient graph relative to $\pi$. Then, for any time $t$,
	\[
		(e^{-itA_X})_{a,b}=(e^{-itB})_{\{a\},\{b\}}
	\]
	and therefore $G$ has perfect state transfer from $a$ to $b$ at time $t$ if and only
    if the symmetrized quotient graph has perfect state transfer from $\{a\}$
	to $\{b\}$.\qed
\end{lemma}

\section{Strongly Cospectral Vertices}

Let $u$ and $v$ be vertices in $X$. We say that $u$ and $v$ are
\textsl{cospectral vertices}
if the characteristic polynomials $\phi(X\diff u,t)$ and $\phi(X\diff v,t)$ are equal.
If $\seq \theta1m$ are the distinct eigenvalues of $X$ and the matrices $\seq E1m$
are the orthogonal projections onto the corresponding eigenspaces we have the spectral
decomposition
\[
	A = \sum_r \theta_r E_r.
\]
From \cite{cg-pst} we know that $u$ and $v$ are cospectral if and only
if $(E_r)_{u,u}=(E_r)_{v,v,}$ for all $r$. Since $(E_r)_{u,u}=\norm{E_re_u}^2$,
we see that $u$ and $v$ are cospectral if and
only if the projections $E_re_u$ and $E_re_v$ have the same length for each $r$. We say
that $u$ and $v$ are \textsl{strongly cospectral} if, for each $r$,
\[
	E_r e_v = \pm E_r e_u.
\]
\cite{Kay-cos} observed that if there is perfect state transfer from vertices
$u$ to $v$, then $u$ and $v$ are strongly cospectral.
In \cite{cg-pst}, an argument due to Dave Witte is presented, which shows that if
there is pretty good state transfer from vertex $u$ to vertex $v$,
then $u$ and $v$ are strongly cospectral. If the eigenvalues of $A$ are simple, then
two vertices are strongly cospectral if and only if they are cospectral.

\begin{lemma}
	\label{lem:cospectral}
	Vertices $u$ and $v$ are strongly cospectral if and only if for each
	idempotent $E_r$, we have $(E_r)_{u,u} = (E_r)_{v,v} = \pm(E_r)_{u,v}$.
\end{lemma}

\proof
The vertices $u$ and $v$ are cospectral if and only if $(E_r)_{u,u} = (E_r)_{v,v}$
for each $r$.
Set $y=E_re_u$ and $z=E_re_v$. By Cauchy-Schwarz
\[
	0 \le \norm{y}^2\norm{z}^2 - |y^Tz|^2
\]
and equality holds if and only if $\{y,z\}$ is linearly dependent.
Since
\[
	\norm{y}^2\norm{z}^2 - |y^Tz|^2 = (E_r)_{u,u}(E_r)_{v,v} -(E_r)_{u,v}^2.
\]
and since $\norm{E_re_u}=\norm{E_re_v}$ when $u$ and $v$ are cospectral, the lemma follows.\qed

We will show that for the symmetric double star graphs $S_{k,k}$, the two
central vertices are strongly cospectral.

We now consider the behaviour of strongly cospectral vertices under
quotients over an equitable partition.

\begin{lemma}
	Let $X$ be a graph and let $\pi$ be an equitable partition of $X$ in which
	$\{u\}$ and $\{v\}$ are singleton cells. Then $u$ and $v$ are strongly cospectral
	in $X$ if and only if $\{u\}$ and $\{v\}$ are strongly cospectral in the
	symmetrized quotient graph $X/\pi$.
\end{lemma}

\proof
Let $Q$ be the normalized characteristic matrix of $\pi$ and let $B$ be the
symmetrized quotient matrix.
If $E_r$ is an idempotent in the spectral decomposition of $A$, then $E_r=p_r(A)$
for some polynomial $p$ and so
\[
    E_rQ = p_r(A)Q = Qp_r(B).
\]
Therefore $p_r(B)=Q^TE_rQ$ and
\[
    AE_rQ = \theta_r E_rQ = \theta_r AQ p_r(B) = \theta_r QBp_r(B),
\]
from which it follows that the non-zero matrices $Q^TE_rQ$ are the idempotents
in the spectral
decomposition of $B$. If $a,b\in\{u,v\}$, then since $\{u\}$ and $\{v\}$
are singleton cells of $\pi$, we have
\[
    (Q^TE_rQ)_{a,b} = (E_r)_{a,b}.
\]
We conclude that $u$ and $v$ are strongly cospectral in $X$ if and only if they are
strongly cospectral in $X/\pi$.\qed

\section{Perfect State Transfer} 

Let  $S_k$ and $S_\ell$ be the two star graphs with $k$ and $\ell$ edges respectively.
Then the double star graph, which we denote by $S_{k,\ell}$, is
the graph obtained by joining the two vertices with degrees $k$ and $\ell$
of $S_k$ and $S_\ell$, respectively.  We call a double star graph
a \textsl{symmetric double star} if $k=\ell$.

In this section we will show that double star graphs do not have perfect state transfer
                                              
\begin{figure}[htbp]
\begin{center}
\includegraphics[totalheight=3cm]{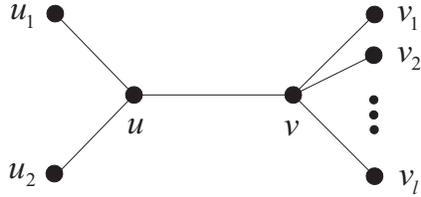}\\
 \caption{\label{path} \small{The double star graph $S_{2,\ell}$}}
\end{center}
\end{figure}

\begin{lemma}\label{lem:Cond-for-S2L}
    Let $S_{2,\ell}$ be a double star graph. Let $\pi$ be the equitable partition with 
    cells 
    \[
        \{\{ \{u_1\}, \{u_2\}, \{u\}, \{v\}, N(v)\setminus \{u\}\}.  
    \]
    If $S_{2, \ell}$ has perfect state transfer between $u_1$ and $u_2$, then 
    $\exp(i \theta_i t) = -1$ for each nonzero eigenvalue $\theta_i$.
\end{lemma}

\proof 
Suppose $u$ and $v$ are the two central vertices of $S_{2,\ell}$, that $u_1, u_2$ are 
the two neighbors of $u$ and  $\seq v1\ell$ are the neighbors of $v$.  Then $\pi$ is an 
equitable partition with cells
\[
	\{ \{u_1\}, \{u_2\}, \{u\}, \{v\}, N(v)\setminus \{u\}\}
\]
and
\[
	B= \left(
	\begin{array}{ccccc}
	0 & 0 & 1& 0 & 0\\
	0 & 0 & 1 & 0 & 0\\
	1 & 1 & 0 & 1 & 0\\
	0 & 0 & 1 & 0 & \sqrt{\ell}\\
	0 & 0 & 0 & \sqrt{\ell} &0
	\end{array}
	\right)
\]
is the adjacency matrix of the corresponding symmetrized quotient graph.
We have perfect state transfer between $u_1$ and $u_2$ if and only if 
the symmetrized quotient graph has perfect state transfer between vertex $\{u_1\}$ 
and $\{u_2\}$. Let
\[
	F=\left(
             \begin{array}{ccccc}
               0 & 1 & 0 & 0 & 0 \\
               1 & 0 & 0 & 0 & 0 \\
               0 & 0 & 1 & 0 & 0 \\
               0 & 0 & 0 & 1 & 0 \\
               0 & 0 & 0 & 0 & 1 \\
             \end{array}
           \right).
\]
Then the symmetrized quotient graph
have perfect state transfer between vertex $\{u_1\}$ and $\{u_2\}$ if and only 
$U_B(t)=\exp(i Bt)= \ga F$, where $\parallel \ga \parallel=1$. The eigenvalues of $B$ are
\begin{align*}
	\th_1&=0,\\
	\th_2&=  \frac{1}{2}\sqrt{2\ell+6+2\sqrt{\ell^2-2\ell+9}},\\[2pt]
	\th_3&=-\frac{1}{2}\sqrt{2\ell+6+2\sqrt{\ell^2-2\ell+9}},\\[2pt]
	\th_4&= \frac{1}{2}\sqrt{2\ell+6-2\sqrt{\ell^2-2\ell+9}},\\[2pt]
	\th_5& = -\frac{1}{2}\sqrt{2\ell+6-2\sqrt{\ell^2-2\ell+9}}.
\end{align*}
Let $E_i$ be the idempotent of $\th_i$, where $i= 1, \cdots, 5$.
The symmetrized quotient graph
has perfect state transfer between vertex $\{u_1\}$ and $\{u_2\}$ if and only if
\begin{equation}\label{eq:F-exp}
	U_B(t)=\sum\limits_{i}\exp(i \th_i t)E_i= \ga F,
\end{equation} 
where $\|\gamma\|=1$.
Note that $(-1,1,0,0,0)$ is an eigenvector for $\th_1$ and
\[
	E_1= \frac{1}{2}\left(
	\begin{array}{ccccc}
	1 & -1 & 0 & 0 & 0\\
	-1 & 1 & 0 & 0 & 0\\
	0 & 0 & 0 & 0 & 0\\
	0 & 0 & 0 & 0 & 0\\
	0 & 0 & 0 & 0 & 0
	\end{array}
	\right)
\]
This implies that $I-2E_1 = F$. Note that $ \sum\limits_{i}E_i=I$ and
\[
	F=-E_1+E_2+E_3+E_4+E_5.
\]
This, together with \eqref{eq:F-exp} and the fact $\exp(i \th_1 t)= 1$, implies that $\ga=-1$ and
$\exp(i \th_i t) = -1$ for $i=2,3,4,5$.\qed

\begin{lemma}\label{lem:no-pst-s2l}
    There is no perfect state transfer between $u_1$ and $u_2$ in the double star graph $S_{2,\ell}$.
\end{lemma}

\proof 
By Lemma \ref{lem:Cond-for-S2L}, $S_{2,\ell}$ has perfect state transfer between 
$u_1$ and $u_2$ if and only if
$ \exp(i \th_i t) = -1 $ for $\th_i \neq 0$, where $\th_i$ is as given in the 
proof of Lemma~\ref{lem:Cond-for-S2L}.
Hence $\th_i t = (2k_i+1)\pi$ for some $k_i \in \ints$.
In particular,
\begin{equation}\label{eq:s2l-RL}
\frac{\th_2}{\th_4}=\frac{2k_2+1}{2k_4+1}
\end{equation}
The right hand of this equation is rational. Next we show that the left hand is irrational.
\begin{align*}
\frac{\th_2} {\th_4} &= \sqrt{\frac{  2\ell+6+2\sqrt{\ell^2-2\ell+9 }}
    {  2\ell+6-2\sqrt{\ell^2-2\ell+9 }  }}\\[2pt]
    &= \frac{\sqrt{2\ell} (l+3+\sqrt{l^2-2l+9}) }{ 4\ell}
\end{align*}

Note first that
\[
	\ell^2-2\ell+9=(\ell-1)^2+8
\]
If $a$ and $b$ are integers and with $a<b$ and $b^2-a^2$ is even, then $b-a$ is even.
As $(a+2k)^2-a^2=4a(a+k)$, we see that $\ell^2-2\ell+9$ is a perfect square if and
only if $\ell=2$. Hence if $\ell \neq2$ then ${\th_2}/{\th_4}$ is irrational
and thus \eqref{eq:s2l-RL} cannot hold.
This means that we cannot have perfect state transfer from $u_1$ to $u_2$ when $\ell\neq2$.

If $\ell = 2$, then $\th_2=2$, $\th_4=1$ and
\[
	\frac{2k_2+1}{2k_4+1} =2,
\]
which is impossible.
Therefore, we cannot have perfect state transfer between $u_1$ to $u_2$ on 
the double star graph $S_{2,\ell}$.\qed

Next, we study the state transfer between two central vertices $u$ and $v$
on the double star graph $S_{k,k}$.

\begin{figure}[htbp]
\begin{center}
\includegraphics[totalheight=3cm]{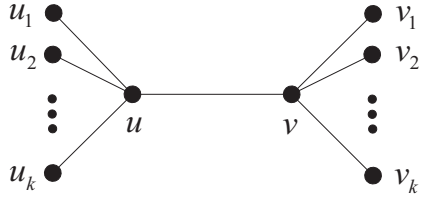}\\
 \caption{\label{path)} \small{The double star graph $S_{k,k}$}}
\end{center}
\end{figure}

\begin{lemma}
	\label{lem:fn23}Let $X$ be the double star graph $S_{k,k}$.
	Then
	\[
		   \exp(i A t)_{u,v}=\bigl((1-2\be)\sin\alpha t+2\be\sin(1-\alpha )t\bigr)i.
	\]
	where $\alpha= \frac{1+\sqrt{1+4k}}{2}$ and $ \be= \frac{k}{1+4k+\sqrt{1+4k}}$.
\end{lemma}

\proof
Let $\pi$ be an equitable partition with cells
$\{N(u)\setminus {v},\{u\},\{v\},N(v)\setminus {u}\}$.
Then the adjacency matrix of the quotient graph $X/\pi$ is 
\[
	A_{X/\pi}= \left (
	\begin {array}{cccc}
	0&1&0&0\\
	k&0&1&0\\
	0&1&0&k\\
	0&0&1&0\\
	\end {array}
	 \right)
\]
After symmetrizing, we get
\[
	B= \left (
	\begin {array}{cccc}
	0&\sqrt{k}&0&0\\
	\sqrt{k}&0&1&0\\
	0&1&0&\sqrt{k}\\
	0&0&\sqrt{k}&0\\
	\end {array}
	 \right)
\]
The eigenvalues of $B$ are $\al, 1-\al, \al-1, -\al$.
The corresponding eigenvectors are the columns of the following
matrix:
\[
\left(
  \begin{array}{cccc}
    -1 & -1 & 1 & 1 \\
   \noalign{\medskip} -\frac{\sqrt{k}}{\al} & \frac{\sqrt{k}}{\al-1}
		& \frac{\sqrt{k}}{\al-1} & -\frac{\sqrt{k}}{\al} \\
  \noalign{\medskip}  \frac{\sqrt{k}}{\al} & -\frac{\sqrt{k}}{\al-1}
		& \frac{\sqrt{k}}{\al-1} & \frac{\sqrt{k}}{\al} \\
  \noalign{\medskip}  1 &1 & 1 & 1 \\
  \end{array}
\right)
\]
Hence
\begin{align*}
	\exp(iBt)_{\{u\},\{v\}} =&(\frac{1}{2}-\be)\exp(i \al t)+\be \exp(i(1-\alpha)t)\\
	&-\be\exp(-i(1-\alpha)t)-(\frac{1}{2}-\be)\exp(-i\alpha t)\\
	&=((1-2\be) \sin(\alpha t)+2\be\sin((1-\alpha )t))i.
\end{align*}
By \lref{gy}, we have
\[
	\exp(iAt)_{u,v}=\exp(iBt)_{\{u\},\{v\}}.
\]
The result follows.\qed

\begin{lemma}\label{lem:nopst}
	Let $S_{k,k}$ be a double star graph. Then there is no perfect
	state transfer from one vertex of degree $k+1$ to the other.
\end{lemma}

\proof
Note that $2\be = \frac{2k}{1+4k+\sqrt{1+4k}}$ and so we have $0<\be<1$. Hence
\begin{align*}
	\parallel\exp(iA t)_{u,v}\parallel &=|(1-2\be)\sin(\alpha t)+2\be\sin((1-\alpha )t)|\\
		&\leq |1-2\be|+|2\be| \\
		&=1.
\end{align*}
Equality holds if and only $\sin\alpha t=\sin(1-\alpha)t=\pm 1$.
Without loss of generality, assume that $\sin\alpha t=\sin(1-\alpha)t=1$. Then
\begin{align}
	\alpha t &= \frac{\pi}{2}+2m\pi \label{eq1},\\
	(1-\alpha)t &=\frac{\pi}{2}+2n\pi.\label{eq2}
\end{align}
It follows that $\alpha=\frac{4m+1}{4(m+n)+2}$. This implies that $\al$ is not an integer.

On the other hand, suppose $\delta=1+4k$ and $\alpha=\frac{1+\sqrt{\delta}}{2}$.
Then we can rewrite
equations (\ref{eq1}) and (\ref{eq2}) in the following form:
\begin{align*}
	\frac{1+\sqrt{\delta}}{2}t &= \frac{\pi}{2}+2m\pi,\\
	\frac{1-\sqrt{\delta}}{2}t &= \frac{\pi}{2}+2n\pi
\end{align*}
and hence
\[
	\frac{1+\sqrt{\delta}}{1-\sqrt{\delta}}=\frac{1+4m}{1+4n}\in \rats.
\]
Note that on the other hand,
\[
	\frac{1+\sqrt{\delta}}{1-\sqrt{\delta}}=\frac{1+\delta+2\sqrt{\delta}}{1-\delta},
\]
which implies that $\delta$ is a perfect square. Since
$\delta=1+4k$ is odd, we can assume that $\delta=(2s+1)^2$, thus
$\alpha=\frac{1+\sqrt{\delta}}{2}=s+1$ is an integer. Contradiction.\qed

Note that if $G$ is a group of automorphisms of the graph $X$ and $u\in V(X)$, then $G_u$
denotes the subgroup consisting of the automorphisms that fix $u$.
We have following result from $\cite{cg-pst}$.

\begin{lemma}
	\label{cddelta}
	Let $G$ be the automorphism group of $X$. If we have perfect state transfer
	from vertex $u$ to vertex $v$ in $X$, then $G_u=G_v$.\qed
\end{lemma}

Next, we will give out main result in this section.

\begin{theorem}
    There is no perfect state transfer on double star graph $S_{k,l}$.
\end{theorem}  

\proof  
Suppose $u$ and $v$ are the two central vertices of $S_{k,\ell}$, that
$\seq u1k$ are the neighbors of $u$ and $\seq v1\ell$ are the neighbors of $v$.
By Lemma \ref{cddelta},
if $k,\ell\ge3$, perfect state transfer can only occur between the
two central vertices. Now if we have perfect state transfer between two
vertices of $X$, they must be cospectral and therefore must have the same degree. 
So $k=\ell$ and perfect state transfer can only occur between two central vertices 
$u$ and $v$. By Lemma~\ref{lem:nopst}, there is no perfect state transfer between 
two central vertices $u$ and $v$ in $S_{k,k}$.

Suppose then that $k$ and $\ell$ are not both at least three, then the graph is 
either the path $P_4$ or $S_{2,\ell}$.
We know that there is no perfect state transfer on $P_4$ in \cite{cddekl-pra}

For $S_{2,\ell}$, we consider two cases. First, 
if $\ell \neq 2$, then by Lemmas~\ref{lem:cospectral} and \ref{cddelta}, 
the only vertices which might have perfect state transfer are vertices with degree one 
that are adjacent to the vertex with degree three. 
By Lemma~\ref{lem:no-pst-s2l}, we know there is no perfect
state transfer in this case.

If $\ell =2$, then perfect state transfer might occur
between two vertices of degree one with a common neighbor, or between
the two vertices of degree three. By Lemma~\ref{lem:no-pst-s2l} and Lemma~\ref{lem:nopst}, 
we also do not have perfect state transfer in this case.\qed

\section{Pretty Good State Transfer}

In this section, we will investigate pretty good state transfer on double
star graph. We say we have pretty good state transfer from $u$ to $v$ if there 
is a sequence $\{t_k\}$ of real numbers and a scalar $\ga$ such that
\[
	\lim_{k\to \infty}U(t_k)=\ga e_v,
\]
where $\norm{\gamma}=1$.

\begin{lemma}
	\label{pgst}
	Let $X$ be a graph and $\pi$ is an equitable partition with $\{u\}$
	and $\{v\}$ are singletons. Then there is pretty good state transfer
	from vertex $u$ to vertex $v$ if and only if there is pretty good state
	transfer from $\{u\}$ to $\{v\}$ in
	the symmetrized quotient graph with adjacency matrix $B$.
\end{lemma}

\proof
By \lref{gy}, we have $\exp(i B t_k)_{\{u\},\{v\}}=\exp(i A t_k)_{u,v}$. Hence
\[
	\lim_{t\to \infty}\exp(i A t_k)_{u,v}=\gamma
\]
if and only if
\[
	\lim_{t\to \infty}\exp(i B t_k)_{\{u\},\{v\}}=\gamma.\qed 
\]

\begin{theorem}\label{lem:pgst-s2l}
	Let $S_{2,\ell}$ be a double star graph, $u_1$ and $u_2$ are two pendent vertices 
	adjacent to the vertex with degree three. Then $S_{2,\ell}$ has pretty good state 
	transfer from vertex $u_1$ to $u_2$ if and only if $l \neq 2$.
\end{theorem} 

\proof 
If $\ell = 2$, then the eigenvalues of $S_{2,2}$ are integers. So in this case 
$S_{2,2}$ is periodic and since we know perfect state transfer does not occur, 
pretty good state transfer does not occur either.

Suppose then that $\ell \neq 2$. By Lemma \ref{lem:Cond-for-S2L}, we see
that there is pretty good state transfer from $u_1$ to $u_2$ if and only if there 
is a sequence of times $(t_{k})_{k \geq 0}$ such that
\begin{equation*}
	\lim_{k \to \infty} \exp{i\th_i t_k} = -1 \ \ \mbox{for}~ i = 2,3,4,5.
\end{equation*}
Note that if $\lim\limits_{k \to \infty} \exp(i \th_2 t_k)=-1$, 
then $\lim\limits_{k \to \infty} \exp(i \th_3 t_k)=-1$.
Similarly, if  $\lim\limits_{k\to\infty} \exp(i \th_4 t_k)=-1$, 
then $\lim\limits_{k \to \infty} \exp(i \th_5 t_k)=-1$.
The question becomes whether we can chose integers $m,n$ such that
\begin{equation*}
	\frac{\th_2}{\th_4} \approx \frac{2m+1}{2n+1},
\end{equation*}
that is, whether we can choose integers $m,n$ such that
\[
	m \th_4 +n (-\th_2) \approx \frac{1}{2}(\th_2 - \th_4)
\]
If $\th_2$ and $\th_4$ are linearly independent over $\mathbb{Q}$, then by Kronecker's
approximation theorem the set
\[
	\{m\th_4 -n \th_2: m,n\in \mathbb{Z}\}
\]
is dense in $\re$. Therefore we can chose $m,n\in\mathbb{Z}$, such that
\[
	\frac{\th_2}{\th_4} \approx \frac{2m+1}{2n+1}
\]
This implies we can choose a series $\{t_k\}_{k \geq 0}$ such that $\lim\limits_{k \to \infty} \exp(i \th_i t_k)=-1$ for all $i = 2, 3,4, 5$.

Next, we show that $\th_2$ and $\th_4$ are linearly independent over $\mathbb{Q}$. 
Assume that for some $x,y \in \mathbb{Q}$,
\[
	x \th_2+y\th_4=0
\]
Then
\[
	\frac{\th_2}{\th_4} = -\frac{y}{x}
\]
is rational. However, by the proof of Lemma \ref{lem:no-pst-s2l}, we know 
${\th_1}/{\th_2}$ is irrational. Contradiction.\qed

Next, we will investigate the pretty good state transfer on double star $S_{k,k}$.

\begin{theorem}
	\label{lem:pgst}
	Let $S_{k,k}$ be a symmetric double star graph. Then either
    \begin{enumerate}[(a)]
		\item
		$1+4k$ is not a perfect square  and there is pretty good state transfer
		from vertex $u$ to vertex $v$, or
	    \item
		$1+4k$ is a perfect square and there is not pretty good state transfer
		from vertex $u$ to vertex $v$.
	\end{enumerate}
\end{theorem}

\proof
If $4k+1$ is a perfect square then the eigenvalues of $S_{k,k}$ are rational
and hence integers.
So in this case $S_{k,k}$ is periodic and since we know that perfect state transfer
does not occur, pretty good state transfer does not occur either.

Assume then that $4k+1$ is not a perfect square. By Lemma \ref{lem:fn23} we see that
there is pretty good state transfer from $u$ to $v$ if and only if there is a
sequence of times $(t_\ell)_{\ell\ge0}$ such that
\begin{equation}
	\label{eq:limsin}
	\lim_{\ell\to\infty}\sin \alpha t_\ell= \lim_{\ell\to\infty}\sin (1-\alpha)t_\ell = \pm1.
\end{equation}

Note that if $\lim_{\ell\to\infty}\sin \alpha t_\ell=\pm1$ then
$\lim_{\ell\to\infty}\cos \alpha t_\ell=0$. Since
\[
	\cos t_\ell = \cos\alpha t_\ell\cos(1-\alpha)t_\ell
		-\sin \alpha t_\ell \sin(1-\alpha)t_\ell
\]
we see that if \eqref{eq:limsin} holds then $\lim_{\ell\to\infty}\cos t_\ell=-1$.

This implies that $t_\ell \approx (2m+1)\pi $ and
$\al t_\ell\approx n\pi+\frac{\pi}{2}$ for $m,n\in \ints$.
The question becomes whether we can choose integers $m,n$ such that
$(2m+1)\alpha-n\approx\frac{1}{2}$.

As $1+4k$ is not a perfect square, $\alpha$ is an irrational number. So $\alpha$ and
$\frac{1}{2}$ are linearly independent over the rationals, and hence by Kronecker's
approximation theorem the set
\[
	\{m\alpha -n: m,n\in \mathbb{Z}\}
\]
is dense in $\re$. Therefore we can chose $m,n\in\mathbb{Z}$, such that
\[
	m\alpha-n\approx\frac{1}{4}-\frac{1}{2}\al.
\]
This implies we can choose a series $\{t_\ell\}$ such that both
$\lim_{\ell\to\infty}\sin \alpha t_\ell=\pm 1$
and $\lim_{\ell\to\infty}\cos t_\ell =-1$.\qed

\section{Recurrence}

If $A$ is the adjacency matrix of a graph $X$ then the set
\[
    \{U(t): t\in\re\}
\]
is an abelian group, which we denote by $G$. In fact $G$ is a 1-parameter subgroup of
the unitary group
on $\cx^n$, where $n=|V(X)|$, and so its closure $\comp{G}$ is an abelian Lie group.
Since $G$ is connected, so is its closure and therefore $\comp{G}$ is isomorphic to a
direct product of some number of copies of $\re/\ints$.
Asking whether there is perfect state transfer from $u$ to $v$ is equivalent to asking
whether there is a matrix $M$ in $G$ such that
\begin{equation}
    \label{eq:muumvv}
    M_{u,u} = M_{v,v} = 0,\quad M_{u,v} = M_{v,u} =\ga
\end{equation}
where $|\ga|=1$. Asking whether there is pretty good state transfer is asking whether
there is a matrix $M$ in $\comp{G}$ such that these conditions hold.

If we prove that pretty good state transfer does occur, we have shown that there is a sequence
of times
$t_\ell$ such that $t_\ell\to\infty$ and, for each $\epsilon>0$ there is time $t_\ell$ such that
$U(t_\ell)$ is within $\epsilon$ of a solution to the conditions of \eqref{eq:muumvv}.
However we can say something more concrete, using the following.

\begin{lemma}
    If $\epsilon>0$, there is a time $T$ such that each element of $\comp{G}$ lies
    within $\epsilon$ of an element of $\{U(s+t): 0\le t\le T\}$, for any $s$.
\end{lemma}

\proof
Define
\[
    S = \{M\in\comp{G}: \norm{M-I} < \epsilon\}
\]
Then $S$ is open and its
translates under the action of $G$ cover $\comp{G}$. Since $\comp{G}$ is compact, some
finite set of translates of $S$ cover $\comp{G}$. Hence there is a time $T$ such that
all elements of $\comp{G}$ lie within $\epsilon$ of an element of $\{U(t): 0\le t\le T\}$.
Since $U(t)$ is unitary it follows that each element of $\comp{G}$ lies within $\epsilon$
of an element of
\[
    \{ U(t+s): 0\le t\le T \},
\]
for any $s$. We conclude that any element of $\comp{G}$ lies within $\epsilon$ of
an element of $\{U(t+s): 0\le t\le T\}$, for any $s$.\qed

\begin{corollary}
    If we have pretty good state transfer from $u$ to $v$ in $X$ then for each positive $\epsilon$
    there is a real number $T$ such that, for each $s$, the set $\{ U(t+s): 0\le t\le T \}$
    contains an element within $\epsilon$ of a solution to \eqref{eq:muumvv}.
\end{corollary}

Since $U(0)=I$, the arguments above also yield the conclusion that, if $\epsilon>0$
then there is a time $T$ such that in each real interval of length $T$ there is
a time $t$ such that $\norm{U(t)-I}<\epsilon$. Thus we can say that any graph
is \textsl{approximately periodic}. (We recall that a graph is \textsl{periodic}
if there is a time $T$ such that $U(t)=\gamma I$, for some complex number $\gamma$
with norm 1. Periodic graphs are studied, and characterized, in \cite{cg-period};
their eigenvalues must be square roots of integers, and consequently these graphs
are rare.)

\section*{Acknowledgement}

We thank Dave Witte Morris for some useful discussion and advice
on 1-parameter subgroups of Lie groups.


\begin{thebibliography}{1}

\bibitem{YG}
\textsc{R.~Bachman, E.~Fredette, J.~Fuller, M.~Landry, M. ~Opperman, C.~Tamon
and A.~Tollefson},
\emph{Perfect state transfer on quotient graphs},
Quantum Information and Computation, 12 (2012), pp. 0293--0313.
[arXiv:1108.0339v2 [quant-ph]]

\bibitem{cddekl-pra}
\textsc{M.~Christandl, N.~Datta, T.~Dorlas, A.~Ekert, A.~Kay, and A.J.~Landahl},
\emph{Perfect transfer of arbitrary states in quantum spin networks},
Phys. Rev. A 71 (2005), p.~032312. [arXiv:quant-ph/0411020v2]

\bibitem{cdel-pst}
\textsc{M.~Christandl, N.~Datta,  A.~Ekert and A.J.~Landahl},
\emph{Perfect state transfer in quantum spin networks},
Phys. Rev. Lett. 92 (2004), p.~187902. [arXiv:quant-ph/0309131v2]


\bibitem{cg-pst}
\textsc{C.~Godsil},
\emph{ State transfer on graphs}, Discrete Mathematics 312 (2012), 129--147.
[arXiv:1102.4898v2 [math.CO]]

\bibitem{cg-period}
\textsc{C.~Godsil},
\emph{Periodic Graphs}, Electron. J. Comb. 18(1) (2011). [arXiv:0806.2074v2 [math.CO ]]


\bibitem{cggrbk}
\textsc{C.~Godsil and G.~Royle},
\emph{Algebraic {G}raph {T}heory},  Springer, New York, 2001.

\bibitem{gkss}
\textsc{C.~Godsil, S.~Kirkland, S.~Severini, J.~Smith},
\emph{Number-theoretic nature of communication in
quantum spin chains}, arXiv:1201.4822v1 [quant.ph].

\bibitem{Kay}
\textsc{A.~Kay},
\emph{Perfect, efficient, state transfer and its application as a constructive tool},
International Journal of Quantum Information 8(4) (2010) 641--676. [arXiv: 0903.4274v3 [quant-ph]]


\bibitem{Kay-cos}
\textsc{A.~Kay},
\emph{Basics of Perfect communication through quantum networks},
Phys. Rev. A 84 (2011), p.~022337. [arXiv:1102.2338v2 [quant-ph]]

\bibitem{viv}
\textsc{V.~Kendon},
\emph{Quantum walks on general graphs},
International Journal of Quantum Information 4(5) (2006) 791-805. [arXiv:quant-ph/0306140v3]

\bibitem{sasesh}
\textsc{N.~Saxena, S.~Severini, and I.~Shparlinski},
\emph{Parameters of integral circulant graphs and periodic quantum dynamics},
International Journal of Quantum Information 5(3) (2007) 417--430.
[arXiv:quant-ph/0703236v1]

\bibitem{VinZh}
\textsc{Luc Vinet, Alexei Zhedanov},
\emph{Almost perfect state transfer in quantum spin chains}.
arXiv:1205.4680v2 [quant-ph].



\end{thebibliography}
\end{document}